\documentclass[11pt]{article}
\usepackage{epsfig}
\usepackage{amsmath,amssymb,amsthm,enumerate}
\begin{document}
%

%
\newtheorem{thm}{Theorem}
\newtheorem{defn}[thm]{Definition}
\newtheorem{prop}[thm]{Proposition}
\newtheorem{cor}[thm]{Corollary}
\newtheorem{lemma}[thm]{Lemma}
\numberwithin{equation}{section}

\newcommand{\beqn}{\begin{equation}}
\newcommand{\eeqn}{\end{equation}}
\newcommand{\nn}{\nonumber}

\def\R{{\mathbb R}}
\def\S{{\cal S}}
\def \e {{\epsilon}}
\def\cleq{{~\preccurlyeq~}}
\def\Ce{{ \widetilde{C}_{0e}^\infty( \R^n)}}
\def\Cz{{\widetilde{C}_e^\infty( \R^n)}}
\def\Sn{{|S^{n-1}|}}
\def\U{{\cal U}}
\def\M{{\cal M}}
\def\Z{{\cal Z}}
\def\Lwrn{{L^2_{e,w}(\R^n)}}
\def\Lp{{L^2_+(\R^n)}}
\def\chie{{\chi_\epsilon}}
\def\wchie{{ \widehat{\chie}}}

\def\Rnp{{ \R^{n-1} \times [0, \infty)}}

\def\pa{{\partial}}
\def \dtsq {{ \frac{1}{2t} \frac{\pa}{\pa t}}}
\def \dxnsq {{ \frac{1}{2x_n} \frac{\pa}{\pa x_n} }}
\def\txn {{ t^2 - x_n^2}}

\def \ep {{\epsilon}}
\def \sxn {{\text{sign}(\xi_n)}}
\def \v {{\widehat{v}}}
\def \g {{\widehat{g}}}
\def \Nc{{\cal N}}
\def \Nh {{\widehat{N}}}

\def\L{{L_e(\R^n, |x_n|^{-1} dx)}}
\def\Co{{\cal C}}

\def \ft {{\tilde{f}}}

\def \f {{\hat{f}}}
\def\C{{\mathbb{C}}}
\newcommand{\la}{\langle}
\newcommand{\ra}{\rangle}

\title{ Spherical means with centers on a hyperplane in even dimensions}

\author{
E K Narayanan\\
Department of Mathematics\\
Indian Institute of Science\\
Bangalore, India - 560012\\
~\\
Email: naru@math.iisc.ernet.in\\
\and
Rakesh\\
Department of Mathematical Sciences\\
University of Delaware\\
Newark, DE 19716, USA\\
~\\
Email: rakesh@math.udel.edu }
\date{January 27, 2010}
\maketitle



\begin{abstract}
Given a real valued function on $\R^n$ we study the problem of recovering the function from its spherical means over spheres centered on a hyperplane. An old paper of Bukhgeim and Kardakov derived an inversion formula for the odd $n$ case with great simplicity and economy. We apply their method to derive an inversion formula for the even $n$ case. A feature of our inversion formula, for the even $n$ case, is that it does not require the Fourier transform of the mean values or the use of the Hilbert transform, unlike the previously known inversion formulas for the even $n$ case. Along the way, we extend the isometry identity of Bukhgeim and Kardakov for odd $n$, for solutions of the wave equation, to the even $n$ case.
\end{abstract}

\section{Introduction}

Below $n>1$ will be a positive integer, a point $x \in \R^n$ will occasionally be written as $x=(x', x_n)$ with $x' \in \R^{n-1}$, $x_n \in \R$, and $|S^{n-1}| = \dfrac{ 2 \pi^{n/2}}{\Gamma(n/2)}$ will be the surface area of the unit sphere in $\R^n$. We take
$\hat{f}(\xi) = \int_{\R^n} e^{-i x \cdot \xi} f(x) \, dx$ as the definition of the Fourier transform of $f$. Define $\Cz$ to be the set of smooth functions $f(x)$ on $\R^n$ which
vanish\footnote{ that is $(\pa^\alpha f)(x)=0$ on $x_n=0$ for all multi-indices $\alpha$} to infinite order on $x_n=0$ and are even in $x_n$. Let $\Ce$ denote the subspace of compactly supported functions in $\Cz$. Finally, define $\Lwrn$ to be the Hilbert space of measurable functions $f(x)$ on $\R^n$ which are even in $x_n$ and for which $ \int_{\R^n} \frac{ |f(x)|^2}{|x_n|} \, dx$ is finite - we use the associated inner product.

 For any continuous function $f(x)$ on $\R^n$ define the spherical mean value operator
\beqn
(Mf)(x,t) = \frac{1}{\Sn} \int_{|\theta|=1} f(x+t \theta) \, d \theta,
\qquad (x,t) \in \R^n \times \R;
\label{eq:Mf}
\eeqn
$(Mf)(x,t)$ is the spherical average of $f$ on the sphere centered at $x$, radius $t$. The goal is the recovery of $f$, given $(Mf)(x,t)$ for all $x$ on the hyperplane $x_n=0$ and for all $t \in \R$.
Clearly on $x_n=0$, $(Mf)(x,t)$ does not change if $f(x)$ is replaced by $f(x',-x_n)$, its reflection across $x_n=0$. Hence one can recover at most the even, in $x_n$, part of $f$. If $f \in \Ce$ then $(Mf)(x,t)$ is an even function of $x_n$ and $t$ and $(Mf)(x',0,t)$ vanishes to infinite order on $t=0$. So we define the map
$\M \, : \, \Ce \rightarrow \Cz$ by
\[
(\M f)(x',t) := (Mf)(x',0,t), \qquad (x',t) \in \R^n.
\]
The goal of this article is to construct an inverse of $\M$ and characterize the range of $\M$. We will do this indirectly.

There is a well known relationship between spherical mean values and solutions of the wave equation. For any $f \in \Ce$, let $u(x,t)$ be the solution of the initial value problem
\begin{gather}
u_{tt} - \Delta u = 0, \qquad \text{on} ~ \R^n \times \R,
\label{eq:ude}\\
u(x,0)=f(x),~ u_t(x,0)=0, \qquad x \in \R^n.
\label{eq:uic}
\end{gather}
Note that $u(x,t)$ is even in $t$ and also in $x_n$ because $f(x)$ is even in $x_n$.
From page 682 of [CH62] there is a formula\footnote{there seems to be a misprint in the formula in [CH62] for the even $n$ case} for $u(x,t)$ in terms of $Mf$ and in terms of $f$. For any smooth even function $h(t)$ on $\R$ define the differential operator
\[
(Dh)(t) := \frac{1}{2t} \frac{ \pa h}{\pa t}, \qquad t \in \R.
\]
Note that $Dh(t)$ makes sense for all $t$ because $h$ is even and $(Dh)(t)$ is an even function of $t$. Further, $D$ may be considered as differentiation with respect to $t^2$ because it may be verified that $D (h(t^2)) = h'(t^2)$.
We have the following relations for all $(x,t) \in \R^{n+1}$, $t \neq 0$: for odd $n$
\begin{align}
u(x,t) &= \frac{ \sqrt{\pi}}{\Gamma(n/2)} t D^{(n-1)/2} (t^{n-2} (Mf)(x,t) )
\label{eq:uModd}\\
& = \pi^{(1-n)/2}\, t \, D^{(n-1)/2} \int_{\R^n}
\delta( t^2 - |x-y|^2) \, f(y) \, dy,
\label{eq:uodd}
\end{align}
and for even $n$
\begin{align}
u(x,t) &= \frac{2}{\Gamma(n/2)} t \, D^{n/2}
\int_0^t \frac{r^{n-1}}{\sqrt{t^2-r^2}}\, (Mf)(x,r) \, dr
\label{eq:uMeven} \\
& = \pi^{-n/2} \, t \, D^{n/2} \int_{\R^n} \frac{ H(t^2 - |x-y|^2)}{\sqrt{ t^2 - |x-y|^2}} \, f(y) \; dy
\label{eq:ueven}
\end{align}
where $H(s)$ is the Heaviside function.
Using (\ref{eq:ude}) and that $f$ vanishes to infinite order on $x_n=0$, it may be verified that $u(x',0,t)$ vanishes to infinite order on $t=0$.
Define the map $\U \, : \, \Ce \rightarrow \Cz $ by
\[
(\U f)(x',t) := u(x',0,t).
\]
 We construct the inverse of $\U$ and characterize the range of $\U$ which will give us the inverse of $\M$ and implicitly characterize the range of $\M$.

Define the tempered distribution $K(t) := t^{-1/2} \, H(t)$, $ t \in \R$ and let $K^{(m)}(t)$ denote the $m$-th derivative of $K(t)$.
\begin{thm}\label{thm:U}
Suppose $n>1$ is even.
\vspace{-0.2in}
\begin{enumerate}[(a)]
\item (Isometry) For any $f \in \Ce$ we have the isometry
\beqn
\int_{\R^n} \frac{ |f(x)|^2}{|x_n|} \, dx = \int_{\R^{n}}
\frac{ |(\U f)(x',t)|^2}{|t|} \, dx' \, dt,
\label{eq:uisom}
\eeqn
so $\U$ extends to an isometry on $\Lwrn$.
\item (Inversion Formula) For any $f \in \Ce$, $x \in \R^n$, $x_n>0$ we have
\beqn
f(x) = \frac{(-1)^{\frac{n}{2}} 2 x_n}{\pi^{n/2}} \,
 \int_{\R^{n-1}} \int_0^\infty \frac{ H(t^2-x_n^2 - |y'|^2)}{\sqrt{t^2-x_n^2 - |y'|^2}} \,
t \, \left ( \dtsq \right )^{n/2}
 \frac{(\U f)(x'+y',t)}{t}  \; dt \, dy'.
\label{eq:Uinv}
\eeqn
\item (Range Characterization) The range of $\U$ is a closed subset of $\Lwrn$. Further, if $F \in \Lwrn$ is in the range of the restricted $\U$ (domain is $\Ce$) then
\begin{enumerate}[(i)]
\item $F(x',t) \in \Cz$ and $H(\tau) \,F(x', \sqrt{\tau}) / \sqrt{\tau} \in H^s(\R^n)$ for all real $s$,
\item
$|\pa_{x'}^\alpha \pa_t^j F(x',t)| \leq C_{j,\alpha} (1+|t|)^{-(n-1)/2}$ for all $j \geq 0$ and all multi-indices $\alpha$.
\end{enumerate}
Finally, if $F$ satisfies (i), (ii) then $F$ is in the range of $\U$ iff $\forall x \in \R^n, ~ x_n \neq 0$, we have
\beqn
\int_{\R^{n-1}} \int_0^\infty \frac{ H(t^2+x_n^2 - |y'|^2)}{\sqrt{t^2+x_n^2 - |y'|^2}} \,
t \,  \left ( \dtsq \right )^{n/2}\frac{F(x'+y',t)}{t}  \; dt \, dy' =0
.
\eeqn
\end{enumerate}
\end{thm}
In the statement of Theorem \ref{thm:U} we were careful to write only absolutely convergent integrals
which forced us into a long-winded statement of (c). A briefer version of (b) and (c) would be that for
all $f \in \Lwrn$
\beqn
f(x) = \frac{x_n}{ \pi^{n/2}} \int_{\R^n} K^{(n/2)}(t^2 - x_n^2 - |y'|^2)
\, (\U f)(x'+y',t) \, dy' \, dt,
\qquad \forall x \in \R^n, ~ x_n >0,
\label{eq:infb}
\eeqn
and an $F \in \Lwrn$ is in the range of $\U$ iff
\beqn
\int_{\R^n} K^{(n/2)}(t^2 + x_n^2 - |y'|^2)
\, F(x'+y',t) \, dy' \, dt \; = \; 0,
\qquad \forall x\in \R^n, ~ x_n \neq 0.
\label{eq:infc}
\eeqn
However, now (\ref{eq:infb}) and (\ref{eq:infc}) have to be interpreted in a sense buried 
in the proof of Theorem \ref{thm:U}.
The result for odd $n$ given in [BK78] has a form similar to the formal expressions above. 
They can be
rewritten in terms of absolutely convergent integrals similar to the expressions in Theorem \ref{thm:U}.

From the hypothesis of Theorem \ref{thm:U} it would seem that the
inversion formula is not applicable to an $f \in C_0^\infty(\R^n)$
which is even in $x_n$ but which does not have the right decay
near $x_n=0$. However, for such $f$, if $u(x',0,t)$ is known then
using $u_{x_n}(x',0,t) =0$ (because $u$ is even in $x_n$) 
and (\ref{eq:ude}) one may determine
$f(x',0)=u(x',0,0)$ and all the $x_n$ derivatives of $f(x) =
u(x,0)$ on $x_n=0$. Knowing these derivatives of $f(x)$ on $x_n=0$
we can create a known function $g(x) \in C_0^\infty(\R^n)$ so that
$\tilde{f} =  f-g$ vanishes to a chosen finite order (or infinite
order) on $x_n=0$. Since $g$ is known we can solve (\ref{eq:ude}),
(\ref{eq:uic}) with $g$ replacing $f$ and hence from $u(x',0,t)$
we can obtain the trace on $x_n=0$ of the solution of
(\ref{eq:ude}), (\ref{eq:uic}) with $f$ replaced by $\tilde{f}$.
Now we can apply the inversion formula to the $\tilde{f}$ data and
recover $\tilde{f}$ and hence $f$.

The early work on the inversion of $\U$ (and hence of $\M$) is
described in [Jo55]. For odd $n$, a nice inversion formula for
$\U$ (and hence $\M$) was given by Bukhgeim and Kardakov [BK78].
Later, inversion formulas for $\M$ for odd and even $n$ were
derived in [No80], [Fa85], [An88], [NRT95], [Kl03],  and [Be09];
[SQ05] has a numerical inversion scheme. 
The problem under consideration has applications to
inverse problems in elasticity and other areas - see [BK78] and
[Fa85]. 
The article [FPR04] studies the similar 
problem of recovering a compactly supported function
from its spherical averages over spheres centered on the boundary
of a region containing the support of the function.

Beltukov's formula in [Be09], for the odd $n$ case, is essentially
the inversion formula in [BK78]. Beltukov rediscovered some of the
ideas in [BK78] (he seems not to have been aware of [BK78]) but
not all the ideas in [BK78]. Applying the ideas in [BK78] to the
even $n$ case, we give a concise derivation of an inversion
formula for the even $n$ case which does not require the Fourier
transform of the data or the use of the Hilbert transform. In
[BK78] an isometry identity was established for solutions of the
wave equation for the odd $n$ case; we extend that identity to the
even $n$ case. 

The proof of Theorem \ref{thm:U} uses the tempered distribution
$N(x) := \pi^{-n/2} K^{(n/2)}(x_n - |x'|^2)$ and its Fourier transform.
\begin{prop}\label{prop:ft}
For even $n>1$ we have
$
\widehat{N}(\xi)  = e^{ i |\xi'|^2/(4 \xi_n) }
$
for all $\xi \in \R^n$.
\end{prop}

To prove (c) in Theorem \ref{thm:U} and to be sure that the integral in (b) is absolutely convergent
we need some decay estimates for solutions of the wave equation - see Theorem 1.1 in [So08].
\begin{prop}\label{prop:decay}
Suppose $n>1$ is a positive integer, $f,g \in C_0^\infty(\R^n)$
and $v(x,t)$ the solution of
\begin{gather}
v_{tt} - \Delta v = 0, \qquad \text{on} ~ \R^n \times \R,
\label{eq:vde}\\
v(x,0)=f(x),~ v_t(x,0)=g(x), \qquad x \in \R^n
\label{eq:vic}
\end{gather}
then
\begin{align}
|v(x,t)| & \cleq (1+|t|)^{-(n-1)/2}, \qquad \text{$n$ odd}
\nn
\\
|v(x,t)| & \cleq (1+|t|)^{-(n-1)/2} \, (1 + ||t|-|x||)^{-(n-1)/2},
\qquad \text{$n$ even}
\label{eq:vdecay}
\end{align}
with the constant determined by $f,g,n$.
\end{prop}
For $n$ odd, $v(x,t)$ is supported in $|t-|x|| \leq R$ if $f,g$ are supported in $|x| \leq R$. This is not true for $n$ even, instead we have an extra decay contribution over the region $|t-|x|| \geq R$.
We have included the statement and the proof for the odd $n$ case because of its possible application to the inversion formula for odd $n$ in [BK78].

The rest of the article is as follows: Section \ref{sec:thm} contains the proof of Theorem \ref{thm:U}, Section \ref{sec:ft} contains the proof of Proposition \ref{prop:ft}, and Section \ref{sec:decay} contains the proof of Proposition \ref{prop:decay}. We have included the proof of Proposition \ref{prop:decay} for completeness though it is stated as a theorem in [So08] whose proof was left as a good exercise.
The proof of Theorem \ref{thm:U} would be quite short (see [BK78] for the odd $n$ case) if we accept formal arguments and expressions; our proof is somewhat long because of the need to show the absolute convergence of the integrals in the expressions.

This work was done while the first named author was supported in part by a grant from the UGC via DSA-SAP and the second named author was on sabbatical at the Department of Mathematics of the Indian Institute of Science in Bangalore, India. The second named author thanks the department for its warm hospitality and generous financial support and the University of Delaware for granting the sabbatical.

\section{Proof of Theorem \ref{thm:U}}\label{sec:thm}
For $t>0$ and $x' \in \R^{n-1}$ from (\ref{eq:ueven}) we have (integrals below to be interpreted as distributional action)
\begin{align*}
\frac{(\U f)(x',0,t)}{t} = \frac{u(x',0,t)}{t} & = \pi^{-n/2} \, \int_{\R^n} K^{(n/2)}(t^2 - |x'-y'|^2 - y_n^2) \; f(y) \; dy
\\
&= 2 \pi^{-n/2} \int_{\R^{n-1}} \int_0^\infty K^{(n/2)}(t^2 - |x'-y'|^2 - y_n^2) \, f(y) \, dy_n \, dy'
\\
&= \pi^{-n/2} \int_{\R^{n-1}} \int_0^\infty K^{(n/2)}(t^2 - |x'-q'|^2 - q_n) \, \frac{f(q', \sqrt{q_n})}{\sqrt{q_n}} \, dq_n \, dq'.
\end{align*}
Hence taking $x'=p'$, $t= \sqrt{p_n}$ for any $p' \in \R^{n-1}$, $p_n >0$, we have
\beqn
\frac{ (\U f)(p',\sqrt{p_n})}{\sqrt{p_n}}
= \pi^{-n/2} \int_{\R^{n-1}} \int_0^\infty K^{(n/2)}(p_n - |p'-q'|^2 - q_n) \, \frac{f(q', \sqrt{q_n})}{\sqrt{q_n}} \, dq_n \, dq'.
\label{eq:ufsq}
\eeqn

Let $\Lp$ denote the subspace of functions $g(x)$ in $L^2(\R^n)$ which are zero for $x_n \leq 0$. Define the specialized zero conversion operator $\Z \; : \; \Lwrn \rightarrow \Lp$ with
\[
(\Z f)(p) =
\begin{cases}
\frac{  f(p',\sqrt{p_n})}{\sqrt{p_n}} & p_n >0,\\
0 & p_n \leq 0.
\end{cases}
\]
Then $\Z$ is an isomorphism from $\Lwrn$ to $\Lp$ and its inverse is given by $(\Z^{-1} g)(x',x_n) = x_n^2 g(x', x_n^2)$. Further $\Z$ maps $\Ce$ to the subspace of functions in $\C_0^\infty(\R^n)$ which are zero on $x_n \leq 0$.

For any $f \in \Ce$ let $g=\Z f$; then from (\ref{eq:ufsq}) we have
\beqn
(\Z \U f)(p) = \pi^{-n/2} \int_{\R^{n}} K^{(n/2)}(p_n - |p'-q'|^2 - q_n) \, g(q) \, dq
= (N*g)(p), \qquad p \in \R^n;
\eeqn
note that this is true even if $p_n \leq 0$.
This suggests we define a map $\Nc \, : \, C_0^\infty( \R^n)
\rightarrow C^\infty( \R^n)$ with
\[
(\Nc g)(p) = (N*g)(p), \qquad p \in \R^n;
\]
then $\Z \U f = \Nc \Z f$. Using
Proposition \ref{prop:ft} and the definition of $\Nc$ we have
\begin{align}
\widehat{\Nc g}(\xi)  =  \widehat{N}(\xi) \, \widehat{g}(\xi)
= e^{ i |\xi'|^2/(4 \xi_n) } \, \widehat{g}(\xi).
\label{eq:vhat}
\end{align}
So $|(\widehat {\Nc g})(\xi)|
= |\widehat{g}(\xi)|$ for all $\xi \in \R^n$ and $\Nc$ extends to an isometric isomorphism from $L^2(\R^n)$ to $L^2(\R^n)$.
In particular, for all $f \in \Ce$ we have $\|\Z f\|_{L^2(\R^n)} = \|\Z \U f\|_{L^2(\R^n)}$ and hence
\[
\int_{\R^{n-1}} \int_0^\infty \frac{| f(p', \sqrt{p_n})|^2}{p_n} \, dp_n \, dp' = \int_{\R^{n-1}} \int_0^\infty \frac{| (\U f)(q', \sqrt{q_n})|^2}{q_n} \, dq_n \, dq'.
\]
Introducing the change of variables $p'=x', p_n = x_n^2$, and $y'=q', t = \sqrt{q_n}$ on the regions $p_n>0$ and $q_n>0$ we obtain (\ref{eq:uisom}) which proves (a) of Theorem \ref{thm:U}.

The equation (\ref{eq:vhat}) allows us to express any $g \in L^2(\R^n)$ in terms of $\Nc g$ giving us
\beqn
\g(\xi) = e^{ -i |\xi'|^2/(4 \xi_n) } \, \widehat{\Nc g}(\xi)
= \Nh(\xi', -\xi_n) \, \widehat{\Nc g}(\xi).
\label{eq:gtemp}
\eeqn
Now, for any $g \in C_0^\infty(\R^n)$, $\g(\xi)$ is a rapidly decaying function, and $\Nh(\xi', \xi_n)$ is bounded, so $ \widehat{\Nc g}(\xi)$ is a rapidly decaying function and hence $\Nc g \in H^m(\R^n)$ for all real $m$. Further $\Nh(\xi', -\xi_n)$ is bounded so $N(q', -q_n) \in H^{-m}(\R^n)$ for any real $m >n/2$; pick one $m$. Then for any $p \in \R^n$, using the Fourier inversion formula and (\ref{eq:gtemp}), we have (below $T_p$ is the translation operator)
\begin{align}
g(p) &= \frac{1}{(2 \pi)^n} \int_{\R^n} \g(\xi) e^{ip \cdot \xi} \, d \xi
= \frac{1}{ (2 \pi)^n} \int_{\R^n} \Nh(\xi', -\xi_n) \, \widehat{\Nc g}(\xi)\,  e^{i p \cdot \xi} \, d \xi
\nn \\
&= \frac{1}{ (2 \pi)^n} \int_{\R^n} \overline{\Nh(-\xi', \xi_n) } \, \,\widehat{T_p \Nc g}(\xi) \, d \xi
= \, \la N(-q', q_n),   (\Nc g)(p+q) \ra
\nn \\
&= \la N(p'-q', q_n -p_n),   (\Nc g)(q) \ra
\label{eq:g2temp}
\end{align}
where the $\la, \ra$ is to be understood as the action of an element in $H^{-m}(\R^n)$ on an element of $H^m(\R^n)$.

For any $f \in \Ce$, take $g=\Z f$, then $\Nc g = \Z \U f$; hence for any $p \in \R^n$, from (\ref{eq:g2temp}) we have
\beqn
(\Z f)(p) = \la N(p'-q', q_n-p_n), (\Z \U f)(q) \ra.
\label{eq:Zftemp}
\eeqn
Working formally, for the moment, (\ref{eq:Zftemp}) implies
\begin{align*}
(\Z f)(p) &= \la K^{(n/2)}(p'-q', q_n - p_n), (\Z \U f)(q) \ra
\\
& = (-1)^{n/2} \la K(p'-q', q_n - p_n),
\left ( \frac{\pa}{\pa q_n} \right )^{n/2} (\Z \U f)(q) \ra
\\
&=
(-\pi)^{n/2} \int_{\R^n} \frac{H(q_n - p_n -|q'-p'|^2)}{ \sqrt{q_n - p_n -|q'-p'|^2}}
\, \left ( \frac{\pa}{\pa q_n} \right )^{n/2} (\Z \U f)( q', q_n) \, dq.
\end{align*}
For any $x \in \R^n$ with $x_n>0$, taking $p'=x', p_n = x_n^2$, and taking $y'=q', t = \sqrt{q_n}$ (note $q_n>0$ on the support of $\Z \U f$), we have
\begin{align*}
\frac{f(x)}{x_n} &= 2 \, (-\pi)^{-n/2}
 \int_{\R^{n-1}} \int_0^\infty \frac{ H(t^2-x_n^2 - |x'-y'|^2)}{\sqrt{t^2-x_n^2 - |x'-y'|^2}} \,
t \,  \left ( \dtsq \right )^{n/2}\frac{(\U f)(y',t)}{t}  \; dt \, dy'
\end{align*}
thus proving (b) of Theorem \ref{thm:U}.

We now address (c) of Theorem \ref{thm:U}. If $f \in \Ce$ then $(\U f)(x',t) \in \Cz$ and from (\ref{eq:vhat}) $\Z \U f \in H^s(\R^n)$ for all real $s$. Further, if $u$ is the solution of (\ref{eq:ude}), (\ref{eq:uic}) then $\pa_x^\alpha \pa_t^m u$ is a solution of (\ref{eq:vde}), (\ref{eq:vic}) for some $f,g \in C_0^\infty(\R^n)$ and hence from Proposition \ref{prop:decay}
\[
|\pa_{x'}^\alpha \pa_t^m (\U f)(x',t)| \leq C_{m,\alpha} (1+|t|)^{-(n-1)/2}
\]
for all $m \geq 0$ and all multi-indices $\alpha$,
proving (i), (ii) of (c) in Theorem \ref{thm:U}.

The extended $\U$ is an isometry and hence its range is closed. A function $F(x',t) \in \Lwrn$ is in the range of $\U$ iff $F= \U f$ for some $f \in \Lwrn$, that is (using the isomorphism of $\Z$) iff $\Z F = \Z \U f =  \Nc \Z f$, that is (using the isomorphism of $\Nc$) iff $\Z f = \Nc^{-1} \Z F$ for some $f \in \Lwrn$. So a necessary condition for $F$ to be in the range of $\U$ is that $(\Nc^{-1} \Z F)(p)=0$ for $p_n < 0$. Conversely,
if $ (\Nc^{-1} \Z F)(p) = 0$ on the region $p_n <0$, then
$g := \Nc^{-1} \Z F$ is in $\Lp$. Let $f := \Z^{-1}g$; then $f \in \Lwrn$ and $g=\Z f$ and hence $\Z f = \Nc^{-1} \Z F$ implying
$\Z F = \Nc \Z f = \Z \U f$, so $F = \U f$ proving our claim.

If $F$ satisfies (i), (ii) of (c) in Theorem \ref{thm:U}, then $\Z F \in H^s(\R^n)$ for all $s$. From (\ref{eq:gtemp}),
$ \widehat{(\Nc^{-1} \Z F)}(\xi) = N(\xi', -\xi_n) \widehat{\Z F}(\xi)$ and since the right side is an integrable function so is the left side, and hence by the Fourier inversion formula, as done for (\ref{eq:g2temp}), we have
\[
(\Nc^{-1} \Z F)(p) = \la N(p'-q', q_n -p_n),  (\Z F)(q) \ra, \qquad p \in \R^n.
\]
Hence, working formally for the moment, $F$ is in the range only if
for all $p' \in \R^{n-1}$, with $p_n>0$ we have (as done earlier)
\begin{align*}
0 &= (\Nc^{-1} \Z F)(p', -p_n) =  \la N(p'-q', q_n +p_n),  (\Z F)(q) \ra
\\
&=  (-\pi)^{n/2} \int_{\R^n} \frac{H(q_n + p_n -|q'-p'|^2)}{ \sqrt{q_n + p_n -|q'-p'|^2}}
\, \left ( \frac{\pa}{\pa q_n} \right )^{n/2} (\Z F)( q', q_n) \, dq
\\
&=(-\pi)^{n/2}
\int_{\R^{n-1}} \int_0^\infty \frac{ H(q_n+p_n - |p'-q'|^2)}{\sqrt{q_n+p_n - |p'-q'|^2}} \,
\left ( \frac{\pa}{\pa q_n} \right )^{n/2} \frac{F(q',\sqrt{q_n})}{\sqrt{q_n}} \; dq_n \; dq' .
\end{align*}
For any $x \in \R^n$ with $x_n \neq 0$, take $p'=x'$, $p_n= x_n^2$. Also, inside the integral use the change of variables $y'=q'$, $t^2 =q_n$. Then $F$ is in the range of $\U$ iff for all $x \in \R^n$ with $x_n \neq 0$ we have
\begin{align*}
0 =
 \int_{\R^{n-1}} \int_0^\infty \frac{ H(t^2+x_n^2 - |x'-y'|^2)}{\sqrt{t^2+x_n^2 - |x'-y'|^2}} \,
t \,  \left ( \dtsq \right )^{n/2}\frac{F(y',t)}{t}  \; dt \, dy'
\end{align*}
proving (c) of Theorem \ref{thm:U}.

So to complete a rigorous proof of the Theorem, for both parts (b) and (c), we have to show that if $F$ satisfies (i), (ii) of part (c) of Theorem \ref{thm:U}, then for all $p \in \R^n$, we have
\begin{align*}
\la N(p'-q', q_n-p_n), (\Z F)(q) \ra
=
(-\pi)^{n/2} \int_{\R^n} \frac{H(q_n - p_n -|q'|^2)}{ \sqrt{q_n - p_n -|q'|^2}}
\, \left ( \frac{\pa}{\pa q_n} \right )^{n/2} (\Z \U f)( q'+p', q_n) \, dq,
\end{align*}
where the inner product $\la \cdot, \cdot \ra$ is to be interpreted as the action of an element in $H^{-m}(\R^n)$ on an element of $H^m(\R^n)$. We now give a proof of this claim.

\begin{figure}[h]
\begin{center}
\epsfig{file=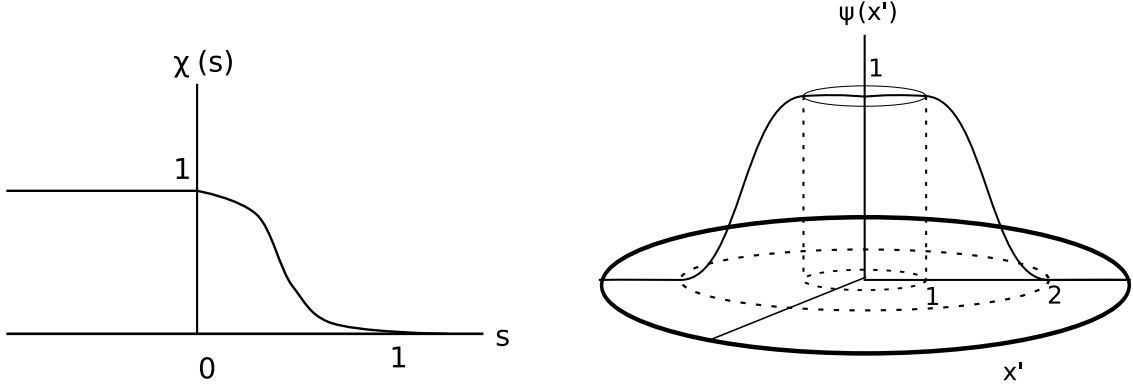, height=2.0in}
\caption{Graphs of $\phi(s)$ and $\psi(x')$}
\label{fig:graphs}
\end{center}
\end{figure}
For any $\phi \in C_0^\infty(\R^n)$, the value of
$\la N(p'-q', q_n-p_n), \phi(q) \ra$, as determined by the duality between $H^{-m}$ and $H^m$, is the same as the value of the action of the tempered distribution  $N(p'-q', q_n-p_n)$ on $\phi(q)$.
Let $\chi \in C^\infty(\R)$ with $\chi(s)=1$ for $s \leq 0$ and $\chi(s)=0$ for $s \geq 1$ and $\psi \in C^\infty(\R^{n-1})$ with $\psi(x')=1$ if $|x'| \leq 1$ and $\psi(x')=0$ if $|x'| \geq 2$.
Then
\[
\phi_k(q) := \chi( \log (q_n^2) -k) \, \psi(q'/k) (\Z F)(q)
\]
 is in $C_0^\infty(\R^n)$ and $\phi_k(q)$ converges to $(\Z F)(q)$ in $H^m(\R^n)$ for all $m$. Hence
\begin{align}
\la N & (p'-q', q_n-p_n),  (\Z F)(q) \ra =
\lim_{k \rightarrow \infty} \la N(p'-q', q_n-p_n), \phi_k(q) \ra
\nn \\
&= \lim_{k \rightarrow \infty}
(-\pi)^{n/2} \la K(q_n - p_n - |p'-q'|^2), (\pa/\pa q_n)^{n/2} \phi_k(q) \ra
\nn \\
&= \lim_{k \rightarrow \infty}
(-\pi)^{n/2} \sum_{j=0}^{n/2} \binom{n/2}{j}
\la K(q_n - p_n - |p'-q'|^2), \pa_n^j(\chi(\log( q_n^2)-k))
\, \psi(q'/k) \, \pa_n^{n/2-j} (\Z F) (q) \ra
\nn \\
&= \lim_{k \rightarrow \infty} (-\pi)^{n/2} \la K(q_n - p_n - |p'-q'|^2), \chi(\log( q_n^2) -k)\, \psi(q'/k) \, \pa_n^{n/2}(\Z F)(q) \ra
\label{eq:NZ2} \\
&= (-\pi)^{n/2} \int_{\R^n} \frac{H(q_n - p_n -|q'|^2)}{ \sqrt{q_n - p_n -|q'|^2}}
\, \left ( \frac{\pa}{\pa q_n} \right )^{n/2} (\Z F)( q'+p', q_n) \, dq.
\label{eq:NZ1}
\end{align}
The equations above required the observation that
\beqn
\lim_{k \rightarrow \infty} \la K(q_n - p_n - |p'-q'|^2),
\pa_n^j (\chi( \log(q_n^2)-k))\, \psi(q'/k) \; \pa_n^{n/2-j} (\Z F) (q) \ra  = 0,
\qquad j>0
\label{eq:obser1}
\eeqn
and that (\ref{eq:NZ2}) equals (\ref{eq:NZ1}). Both these observations follow from the use of the dominated convergence theorem,
that for $j>0$, $\pa_n^j (\chi( \log(q_n^2)-k))$ is supported in $ e^{k/2} \leq |q_n| \leq e^{(k+1)/2}$, and
the integrability in $q$ of $ K(q_n - p_n - |p'-q'|^2) \,
(1+|q_n|)^{-j} \, \pa_n^{n/2  -j} (\Z F) (q)$ for $j=0 \cdots n/2$.

We now prove the integrability assertion in the previous sentence.
For any function $h(t)$ and any real number $k$ we have
$\dtsq (t^k h(t)) = c_1 t^{k-2} h(t) + c_2 t^{k-1} \pa_t h (t)$, so for any $m \geq 0$ we have
\beqn
\left ( \dtsq \right )^m ( t^k h(t) )
= \sum_{l=0}^m c_l t^{k-l - 2(m-l)} \pa_t^l h(t)
=\sum_{l=0}^m c_l t^{k - 2m+l} \pa_t^l h(t).
\label{eq:ph}
\eeqn
Hence using the substitution $q_n=t^2$ we have
\begin{align*}
\left ( \frac{\pa}{\pa q_n} \right )^m (\Z F )(q)
& = \left ( \frac{\pa}{\pa q_n} \right )^m \frac{ F(q', \sqrt{q_n} )}{\sqrt{q_n}}
= \left ( \dtsq \right )^m \frac{ F(q',t)}{t}
= \sum_{l=0}^m t^{-1-2m+l} \pa_t^l F(q',t)
\end{align*}
and from the decay properties of $F$ we conclude that
\beqn
\left | \left ( \frac{\pa}{\pa q_n} \right )^m (\Z F)(q)
\right | \cleq (1+|t|)^{-1-m} (1+|t|)^{-(n-1)/2}|_{t= \sqrt{q_n}}
\cleq (1+|q_n|)^{ - (2m+n+1)/4}.
\label{eq:qdecay}
\eeqn
Hence for any $j \geq 0$,
\begin{align*}
\int_{\R^n} & \frac{K(q_n-p_n - |p'-q'|^2)}{(1+|q_n|)^j}
\, | \pa_n^{(n/2)-j} (\Z F)(q) | \, dq
\\
&~~~~= \int_{p_n}^\infty \int_{|q'|^2 \leq q_n - p_n}
\frac{ | \pa_n^{(n/2)-j} (\Z F)(q'+p', q_n) |}{ \sqrt{q_n - p_n - |q'|^2}} \, (1+|q_n|)^{-j} \, dq' \, dq_n
\\
&~~~~= \int_{p_n}^\infty \int_{|z'|^2 \leq 1}
 \frac{ | \left ( \pa_n^{(n/2)-j}  (\Z F) \right )
 (\sqrt{q_n-p_n} \, z'+p', q_n)  |}{ \sqrt{1 - |z'|^2}} \, (1+|q_n|)^{-j} \, (q_n - p_n)^{(n-2)/2}
\, dz' \, dq_n
\\
&~~~~ \cleq \int_{p_n}^\infty \int_{|z'|^2 \leq 1} \frac{ (1+|q_n|)^{-5/4}}{ \sqrt{1 - |z'|^2} } \, dz' \, dq_n
\end{align*}
which is finite; we have used (\ref{eq:qdecay}) and noted that for all $j=0 \cdots n/2$
\[
(-j + (n-2)/2) - ( n+1 + 2(n/2-j))/4 = -j/2 - 5/4 \leq  -5/4.
\]

\hfill{\bf QED}

\section{Proof of Proposition \ref{prop:ft}}\label{sec:ft}

Let $\S(\R^n)$ denote the usual space of smooth, rapidly decreasing functions on $\R^n$. Choose  $\chi(s) \in C^\infty(\R)$ with $\chi(s) = e^{-s}$ for $s \geq 0$ and $\chi(s) = 0$ for $s \leq -1$; then $\chie$, defined through
\[
\wchie(x) :=  e^{- \e^2 |x'|^2} \, \chi( \e x_n),
\qquad x \in \R^n,
\]
satisfies $\lim_{\epsilon \rightarrow 0^+} \chie*\phi = \phi$ for all $\phi \in \S(\R^n)$ in the topology of $\S(\R^n)$.

If $P(x) = K(x_n- |x'|^2)$ then $N(x) = \pi^{-n/2} \,\pa_n^{n/2} P(x)$ and hence
$\Nh(\xi) = \pi^{-n/2} ( i \, \xi_n)^{n/2} \widehat{P}(\xi)$. By definition, for all
$\phi \in \S(\R^n)$ we have
\beqn
\la \widehat{P}, \phi \ra  = \lim_{\epsilon \rightarrow 0^+} \la \widehat{P}, \chie*\phi \ra
= \lim_{\epsilon \rightarrow 0^+} \la P, \wchie \widehat{\phi} \ra.
\label{eq:Pe}
\eeqn
Now (all integrals below are absolutely convergent)
\begin{align*}
\la P, \wchie \widehat{\phi} \ra
&=
\int_{|x'|^2 \leq x_n} \frac{ e^{- \e^2 |x'|^2} \, \chi( \e x_n)}{\sqrt{x_n - |x'|^2}} \widehat{\phi}(x) \, dx
=
\int_{\R^{n-1}} \int_{|x'|^2}^\infty \int_{\R^n} \frac{ e^{- \e^2 |x'|^2} \, e^{-\e x_n}}{\sqrt{x_n - |x'|^2}} \; e^{- i x \cdot \xi}
\, \phi(\xi) \, d \xi \, dx_n \, dx'
\\
&=
\int_{\R^n} F_\e(\xi) \, \phi(\xi) \, d \xi
\end{align*}
where
\begin{align*}
F_\e(\xi) &= \int_{\R^{n-1}} \int_{|x'|^2}^\infty  \frac{ e^{- \e^2 |x'|^2} \, e^{-\e x_n}}{\sqrt{x_n - |x'|^2}} \; e^{- i x \cdot \xi}
\, dx_n \, dx'
\\
& =\int_{\R^{n-1}} e^{-(\e^2+\e + i \xi_n)|x'|^2} e^{-i x' \cdot \xi'} \, d x'
\;
\int_0^\infty \frac{ e^{-x_n( \e + i \xi_n)}}{\sqrt{x_n}} \, dx_n
\\
&= \frac{ \pi^{(n-1)/2}} { (\sqrt{\e^2 +\e + i \xi_n})^{n-1} }
\, \exp \left ( \frac{- |\xi'|^2}{ 4 (\e^2 +\e + i \xi_n)} \right )
\, \frac{ \sqrt{\pi}}{ \sqrt{ \e + i\xi_n}}
\end{align*}
 and the branch of square root used is the one obtained by cutting out the negative reals and whose argument lies in $(-\pi/2,\pi/2)$.
Hence
\begin{align*}
\la \Nh, \phi \ra &= \la \pi^{-n/2} \, (i \xi_n)^{n/2} \, \widehat{P}(\xi) , \phi(\xi) \ra
= \la \widehat{P}(\xi) , \pi^{-n/2} \, (i \xi_n)^{n/2} \, \phi(\xi) \ra
\\
&= \lim_{\epsilon \rightarrow 0^+}
\int_{\R^n}\pi^{-n/2} \, (i \xi_n)^{n/2} \, F_\e(\xi) \, \phi(\xi) \, d \xi.
\end{align*}
Now $\pi^{-n/2} \, (i \xi_n)^{n/2} \, F_\e(\xi)$ may be seen to be bounded with a bound independent of $\epsilon$ and its pointwise limit is $\exp(i|\xi'|^2/(4 \xi_n))$ almost everywhere, hence by the dominated convergence theorem
\[
\la \Nh, \phi \ra = \int_{\R^n} \exp(i|\xi'|^2/(4 \xi_n)) \, \phi(\xi) \, d \xi.
\]
\hfill{\bf QED}

\section{Proof of Proposition \ref{prop:decay}}\label{sec:decay}

Because of linearity it is enough to deal only with the two special cases where either $f=0$ or $g=0$. We give the proof only of the case where $g=0$, in which case $v$ is given by the formulas (\ref{eq:uodd}), (\ref{eq:ueven}) for $n$ odd, even respectively.
The case where $f=0$ is dealt with similarly with the formula for the new $v$ being almost the same as (\ref{eq:uodd}), (\ref{eq:ueven}) except the $t$ term is dropped, the $D$ power is reduced by one, and the expressions must be multiplied by $2$.

\subsection{Odd $n$ case}
Let $f \in C_0^\infty(\R^n)$ with $f$ supported in the $0$ centered ball of radius $R$.
Since $(Mf)(x,t) \neq 0$ iff the sphere $|y-x|=t$ intersects the ball $|y| \leq R$, we see that $(Mf)(x,t)$ is supported in the region $ | \, t -|x| \,| \leq R$, hence, using (\ref{eq:uModd}), $u(x,t)$ is supported in the region $ | \, t -|x| \,| \leq R$. Now suppose $t>0$, $p$ a real number and $h \in C_0^\infty(\R^n)$; then
\begin{align}
D (t^p \, (Mh)(x,t)) &= \frac{p}{2} \, t^{p-2} \, (Mh)(x,t) +
\frac{t^{p-1}}{2|S^{n-1}|} \int_{|\theta|=1} \theta \cdot (\nabla h)(x+t \theta) \, d \theta
\nn \\
&= \frac{p}{2} \, t^{p-2} \, (Mh)(x,t) +
\frac{t^{p-2}}{2|S^{n-1}|} \int_{|\theta|=1} (x+t \theta - x) \cdot (\nabla h)(x+t \theta) \, d \theta
\nn\\
&= c_1 \, t^{p-2} \, (Mh)(x,t)  + c_2 t^{p-2} M( y \cdot \nabla h(y))(x,t)
+ c_3 t^{p-2} x \cdot M( \nabla h)(x,t) \\
&= t^{p-2} \sum_{|\alpha| \leq 1}  x^\alpha M(h_ \alpha)(x,t)
\label{eq:DM1}
\end{align}
where $h_\alpha(x)$ are smooth functions which are just made of derivatives of $h$, possibly multiplied with powers of $x$.
Hence, using (\ref{eq:DM1}) $(n-1)/2$ times in (\ref{eq:uModd}), we see that
\begin{align}
 u(x,t) &= c \, t \, D^{(n-1)/2} t^{n-2} (Mf)(x,t)
=  \sum_{|\alpha| \leq (n-1)/2} x^\alpha (Mh_\alpha)(x,t)
\label{eq:uDM1}
\end{align}
where $h_\alpha$ are made up of derivatives of $f$, possibly multiplied with powers of $x$.
Now for any $h$ supported in $|x| \leq R$, the surface area of the intersection of the sphere $|y-x|=t$ with the ball $|y| \leq R$ will be bounded above by a constant independent of $x$ or $t$. Hence for $t \geq 1$
\beqn
|(Mh)(x,t)| \leq \frac{1}{t^{n-1}\, |S^{n-1}|} \int_{|y-x|=t} |h(y)| \, dS_y
\;
\cleq \; (1+t)^{-(n-1)}
\label{eq:Mgest}
\eeqn
and noting that $u$ is supported in $|t - |x|| \leq R$, for $t>1$ we have from (\ref{eq:uDM1}) that
\[
|u(x,t)| \cleq \frac{(1+|x|)^{(n-1)/2}}{ (1+t)^{n-1}}
\cleq (1+|t|)^{-(n-1)/2}.
\]

\subsection{Even $n$ case}
First note that $u(x,t)$ is supported in $ |x| \leq t+R$. We will show the decay rates, as $t \rightarrow \infty$, for the two cases $|x|-R \leq t \leq |x| + 3 R$ and
$|x| + 3R \leq t$ separately.

Suppose $|x| + 3R \leq t$; then for any $y$ in the support of $f$ we have $|x-y| \leq |x| + |y| \leq |x|+R \leq t$. Hence the region $|x-y| \leq t$ includes the support of $f$. So from (\ref{eq:ueven}) we have
\begin{align*}
 |u(x,t)| &= \pm c t D^{n/2} \int_{|y| \leq R} \frac{ f(y)}{\sqrt{t^2-|x-y|^2}} \, dy
= \pm c \, t \int_{|y| \leq R} \frac{ f(y)}{(t^2 - |x-y|^2)^{(n+1)/2}} \, dy
\\
&
\cleq  \frac{t}{ (t^2 - (|x|+R)^2)^{(n+1)/2}}
 = \frac{t}{ (t + |x| + R)^{(n+1)/2}} \, \frac{1}{ (t - |x| - R)^{(n+1)/2}}.
\end{align*}
Now for $|x|+3R \leq t$ we have $R \leq (t-|x|)/3$ and
\[
t-|x|-R \geq t-|x| - 2R + R  \geq t-|x| - \frac{2}{3} (t-|x) + R
\geq c (|t-|x||+1),
\]
hence
\[
|u(x,t)| \cleq (1+t)^{-(n-1)/2} \, (1+|t-|x||)^{-(n+1)/2}.
\]
This estimate is better than the estimate in the statement of the proposition - the weaker estimate in the statement of the proposition is tight when $f=0$ and $g$ is non-zero.

For $|x|-R \leq t \leq |x|+3R$ we use a technique similar to one used in the odd case. Since we are interested in situation when $t$ is large, we may assume that $|x| > \max(1,R)$.
From (\ref{eq:ueven}) we have
\begin{align}
u(x,t) &= c t \, D^{n/2}
\int_0^t \frac{r^{n-1}}{\sqrt{t^2-r^2}}\, (Mf)(x,r) \, dr
= ct D^{n/2} t^{n-1} \int_0^1 \frac{s^{n-1}}{\sqrt{1-s^2}}\, (Mf)(x,ts) \, ds.
\label{eq:nMeven}
\end{align}
Now, for $t>0$, for any real number $p$ and any smooth function $h(x)$ we have
\begin{align*}
D \left ( t^p (Mh)(x,st) \right ) &= \frac{p}{2} \, t^{p-2} \, (Mh)(x,st) +
\frac{t^{p-1}}{2|S^{n-1}|} \int_{|\theta|=1} s \theta \cdot (\nabla h)(x+ st \theta) \, d \theta\\
&= \frac{p}{2} \, t^{p-2} \, (Mh)(x,st) +
\frac{t^{p-2}}{2|S^{n-1}|} \int_{|\theta|=1} (x+st \theta - x) \cdot (\nabla h)(x+st \theta) \, d \theta \\
&= c_1 \, t^{p-2} \, (Mh)(x,st)  + c_2 t^{p-2} M( y \cdot \nabla h(y))(x,st)
+ c_3 t^{p-2} x \cdot M( \nabla h)(x,st) \\
&= t^{p-2} \sum_{|\alpha| \leq 1} x^\alpha M(h_ \alpha)(x,st)
\end{align*}
where $h_\alpha(x)$ are smooth functions which are just made of derivatives of $h$, possibly multiplied with powers of $x$.
Hence
\begin{align*}
t D^{n/2} t^{n-1} (Mf)(x,ts)
&=\sum_{|\alpha| \leq n/2} x^\alpha (Mh_\alpha)(x,st);
\end{align*}
using this in (\ref{eq:nMeven}), noting that $(Mh_{\alpha})(x,st)$ is supported in $|x|-R \leq st \leq |x|+R$, and using (\ref{eq:Mgest}) with $h$ replaced by $h_\alpha$, we obtain (note $|x| \geq R$)
\begin{align*}
|u(x,t)| & \cleq |x|^{n/2} \int_0^1 \frac{s^{n-1}}{\sqrt{1-s^2}}
\, |(Mh_\alpha)(x,st)| \, ds
\cleq \frac{|x|^{n/2}}{t^{n-1}} \int_{(|x|-R)/t}^1 \frac{1}{\sqrt{1-s^2}}
\, ds
\\
& = \frac{|x|^{n/2}}{t^{n-1}} \; \arccos ((|x|-R)/t).
\end{align*}
We are considering the region where $t-3R \leq |x| \leq t+R$, and there is no loss of generality in assuming that $t> \max(1, 4R)$; hence
\[
|u(x,t)| \cleq  (1+t)^{-(n-2)/2} \, \arccos ((t-4R)/t)
\cleq (1+t)^{-(n-1)/2}
\]
where the last step follows from
$
\lim_{t \rightarrow \infty} t^{1/2} \, \arccos ((t-4R)/t)
= c \neq 0
$
by L'Hopital's rule. So Proposition \ref{prop:decay} follows from the fact that $ R \leq |t-|x|| \leq 3R$ in the region of interest. 
\hfill{\bf QED}

\end{document}